\documentstyle[amssymb,12pt]{article}
\title{The consistency strength of
projective uniformization, revisited}
\author{Ralf Schindler\footnote{The author gratefully
acknowledges a DFG fellowship. He is heavily indebted to
Kai Hauser, John Steel, and, indirectly, Hugh
Woodin.}}
\newtheorem{defn}{Definition}[section]
\newtheorem{thm}[defn]{Theorem}

\begin{document}

\maketitle
Consider the following assumptions, whose conjunction we 
denote by $(RP)$:

\bigskip
(1) every projective set of reals is Lebesgue measurable and has the
property of Baire, and
 
\bigskip
(2) every projective subset of the plane has a projective
uniformization.

\bigskip
Woodin had asked, in \cite{CSPU}, whether $(RP)$ implies Projective
Determinacy. This is not the case, by a recent observation of Steel:

\begin{thm}
(Woodin, Steel) Suppose $V=K$, where $K$ is Steel's core model.
If there are $\kappa_0 < \kappa_1 < ...$ with supremum
$\lambda$ such that for all $n < \omega$ $\kappa_n$ is $< \lambda^+$
strong [i.e., for all $x \in H_{\lambda^+}$ there is $\pi \colon
V \rightarrow M$ with critical point $\kappa_n$ and $M$ transitive
such that $x \in M$] then $(RP)$ holds in a generic extension.
\end{thm}

We here show that this is best possible:

\begin{thm}
If $ZFC \ + \ (RP)$ holds and 
Steel's $K$ exists then $J^K_{\omega_1} \models$
there are infinitely many strong cardinals.
\end{thm}

{\sc Proof.} Suppose not. Let $n < \omega$ be the number of strongs
in $J^K_{\omega_1}$. We work towards a contradiction.

\bigskip
{\it Case 1.} $\omega_1$ is a successor in $K$.

\bigskip
Then Corollary 2.2 of \cite{proj} gives that $J^K_{\omega_1}$ is
(boldface) \boldmath $\Delta$\unboldmath$^1_{n+4}$ (in the codes). But
then we get a projective sequence of distinct reals of length
$\omega_1$, contradicting \cite{CSPU}.

\bigskip
{\it Case 2.} $\omega_1$ is inaccessible in $K$.

\bigskip 
Let $\Phi_m(M)$ denote the following statement, for $m \geq n$:

\bigskip
$M$ is a countable $m$-full mouse, $M \models$ 
there are $\leq m$
many strongs, and for all countable $m$-full $N$, if $M$, $N$ simply
coiterate to $M^*$, $N^*$ with iteration maps $i \colon M \rightarrow
M^*$ and $j \colon N \rightarrow N^*$ such that $M^*$ is an initial
segment of $N^*$ then $i$''$M \subset j$''$N$.

\bigskip
The concept of $m$-fullness was defined in \cite{proj} where we showed
that $\Phi_m(J^K_\kappa)$ holds for all $\kappa \leq$ 
the $(m+1)^{st}$
strong cardinal of $J^K_{\omega_1}$ 
which is either a double successor or an inaccessible in $K$.

It is also shown in \cite{proj} that if $\omega_1$ is
inaccessible in $K$ and there are $\leq m$ strong cardinals in
$J^K_{\omega_1}$ then $\Phi_m(M)$ characterizes (in a $\Pi^1_{m+4}$
way) (cofinally many of) the proper 
initial segments of $J^K_{\omega_1}$. (Cf. \cite{proj} Theorem 2.1. 
This gives a (lightface) $\Delta^1_{m+5}$ definition of
$J^K_{\omega_1}$.)

In particular, for all $m \geq n$ the following holds, abbreviated by
$\Psi^m_n$:

\bigskip
For any two $M$, $M'$, if $\Phi_m(M)$ and 
$\Phi_m(M')$ both hold then $M$
and $M'$ are lined up and if ${\tilde M}$ is the "union"
of all $M$'s with $\Phi_m(M)$ then $On \cap {\tilde M} =
\omega_1$ and ${\tilde M} \models$ there are exactly 
$n$ strong cardinals.

\bigskip
Notice that $\Psi^m_n$ is $\Pi^1_{m+5}$.

\bigskip
By \cite{CSPU}, there is a model $P = L_{\omega_1}[X] \models ZFC$ for
some $X \subset \omega_1$, such that 

\bigskip
(a) $P[g]$ is $\Sigma^1_{n+1000}$ correct in $P[g][h]$ whenever $g$ is
set-generic over $P$ and $h$ is set-generic over $P[g]$, and

\bigskip
(b) $P$ is $\Sigma^1_{n+1000}$ correct in $V$.

\bigskip
Now $P \models \Psi^{n+94}_n$ by (b) and the fact that 
$\Psi^{n+94}$ holds in $V$. Moreover, $P$ is closed under
the dagger operator by (a), so Steel's $K$ exists in $P$, denoted by
$K^P$, and $K^P \models$ there are $> n$ strong cardinals, by
\cite{CSPA} and (a). We may pick $g$ $Col(\mu,\omega)$-generic over
$P$ for some appropriate $\mu$ such that in $P[g]$, $J^K_{\omega_1}
\models$ there are $> n$ strongs. By (a), $\Psi^{n+94}_n$ still holds
in $P[g]$. 

From this we now derive a contradiction, by working in $P[g]$ for the
rest of this proof. So let us assume that (in $V$) $\Psi^{n+94}_n$
holds, Steel's $K$ exists, and
$J^K_{\omega_1} \models$ there are $> n$ strongs.

By $\Psi^{n+96}_n$, there is a J-model ${\tilde M}$ of height
$\omega_1$ such that ${\tilde M} \models ZF^{-} \ + $ there are
exactly $n$ strong cardinals, and $\Phi_{n+96}(M)$ holds
for every proper initial segment $M$ of ${\tilde M}$. 
By \cite{proj}, there is a universal
weasel $W$ end-extending ${\tilde M}$ such that for all countable (in
$V$) $\kappa$ which are cardinals in $W$ and such that
$J^{\tilde M}_\kappa \models$ there are $< n+94$
strong cardinals, $W$ has the definability property at $\kappa$. [This
follows from the fact that cofinally many 
proper initial segments of ${\tilde M}$
are $n+94$ full].

Because $W$ is universal, there is some $\sigma \colon K \rightarrow
W$ given by the coiteration of $K$ with $W$. Let $\kappa$ denote the
$(n+1)^{st}$ strong cardinal of $J^K_{\omega_1}$. 
By a remark above, $J^K_\kappa$
is an initial segment of $W$. But this implies that the critical point
of $\sigma$ is $> \kappa$. [This follows from the fact that if $\mu$
is strong in $J^K_\kappa$, or $\mu = \kappa$,
then $K$ as well as $W$ has the definability
property at $\mu$.] But now, using $\sigma$, ${\tilde M} =
J^W_{\omega_1} \models$ there are at least $n+1$ many strong
cardinals. Contradiction!

\bigskip
\hfill $\square$


\begin{thebibliography}{99}
\bibitem{CSPA} Hauser, Kai, {\it The consistency strength
of projective absoluteness}.
\bibitem{proj} Schindler, Ralf, {\it The projectiveness
of $K \cap HC$ }, handwritten notes.
\bibitem{CMIP} Steel, John, {\it The core model 
iterability problem}.
\bibitem{CSPU} Woodin, Hugh, {\it The
consistency strength of projective
uniformization}.
\end{thebibliography}
\end{document}